# Lattice trees, percolation and super-Brownian motion

Gordon Slade

April 2, 1999

**Abstract**

This paper surveys the results of recent collaborations with Eric Derbez and with Takashi Hara, which show that integrated super-Brownian excursion (ISE) arises as the scaling limit of both lattice trees and the incipient infinite percolation cluster, in high dimensions. A potential extension to oriented percolation is also mentioned.

# 1 Introduction

This paper concerns lattice trees and percolation on $\mathbb{Z}^d$. ¿From the point of view of statistical mechanics, one of the fundamental problems in the study of these models is the construction and analysis of the scaling limit, in which the lattice spacing goes to zero. Control of the scaling limit is closely related to control of the model's critical exponents. General features of the scaling limit are beginning to emerge [?, ?], but much work remains to be done. In particular, there is still no proof of the existence of a single critical exponent for either model in low dimensions.

However, in high dimensions, there has been recent progress for both models. This progress has relied on the fact that the scaling limits in high dimensions turn out to involve integrated super-Brownian excursion (ISE), a close relative of super-Brownian motion (SBM). SBM is a fundamental example of a measure-valued process, a class of objects that has been intensively studied in the probability literature [?, ?, ?].

## 2  SBM and ISE

We will not give a precise mathematical definition of super-Brownian motion here. Our goal in this section is to introduce the key functions associated with ISE that will appear in the results for lattice trees and percolation. A more detailed description of SBM can be found in the article by Cox, Durrett and Perkins in this volume [?], which describes how SBM arises as the scaling limit also for the voter model and the contact process.

SBM can be constructed as an appropriate scaling limit of a critical branching random walk on $\mathbb{Z}^d$, originating from a single initial particle, in the limit as the lattice spacing is shrunk to zero. Such a construction is described in [?], and defines SBM as a remarkable Markov process in which the state at any particular time is a random finite measure on $\mathbb{R}^d$ representing the mass density of particles present at that time. The process dies out in finite time. The entire family tree of SBM is a random finite measure on $\mathbb{R}^d$ and is referred to as the historical process. For dimensions $d \geq 4$, it is almost surely supported on a subset of $\mathbb{R}^d$ having Hausdorff dimension 4 [?, ?].

The mean measure of SBM at time $t$, which represents the mass density at time $t$ averaged over all family trees, is a deterministic measure that is absolutely continuous with respect to Lebesgue measure in all dimensions $d \geq 1$. In fact, it is the probability measure on $\mathbb{R}^d$ with density

$$p_t(x) = \frac{1}{(2\pi t)^{d/2}} e^{-x^2/2t}. \tag{2.1}$$

This density is the transition density for Brownian motion in $\mathbb{R}^d$ to travel from 0 to $x$ in time $t$.

ISE is the random measure on $\mathbb{R}^d$ obtained by conditioning the historical process to be a probability measure on $\mathbb{R}^d$. Alternately, it can be constructed from critical branching random walk on $n^{-1/4}\mathbb{Z}^d$, starting from a single particle and conditional on a fixed size $n$ for the total size of the initial particle's family tree up to extinction, in the limit $n \to \infty$. The law of ISE is a probability measure $\mu_{\text{ISE}}$ on the space $M_1(\mathbb{R}^d)$ consisting of probability measures on $\mathbb{R}^d$ and equipped with the topology of weak convergence. The mean of $\mu_{\text{ISE}}$ is a deterministic probability measure on $\mathbb{R}^d$, which corresponds to averaging over all family trees resulting from the initial particle, under the basic unit mass condition required by ISE. Define

$$a^{(2)}(x, t) = t e^{-t^2/2} p_t(x). \tag{2.2}$$

The mean ISE measure is absolutely continuous with respect to Lebesgue measure in all dimensions $d \geq 1$. Its density with respect to Lebesgue

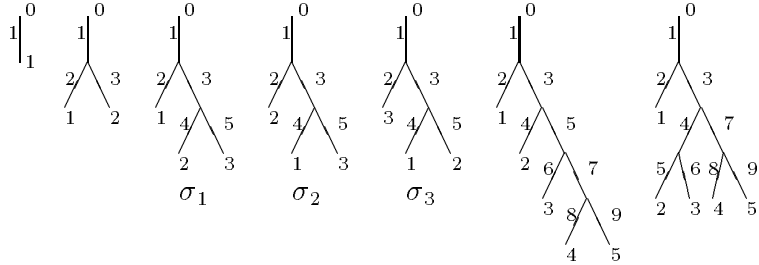

Figure 1: The shapes for $m = 2, 3, 4$, and examples of the $7!! = 7 \cdot 5 \cdot 3 = 105$ shapes for $m = 6$. The shapes' edge labellings are arbitrary but fixed.

measure is the function
$$A^{(2)}(x) = \int_0^\infty a^{(2)}(x, t) dt = \int_0^\infty t e^{-t^2/2} p_t(x) dt. \qquad (2.3)$$

The functions (??) and (??) are ISE two-point functions. A discussion of higher point functions requires the notion of *shape*, which is defined as follows. We start with an *m-skeleton*, which is a tree having $m$ unlabelled external vertices of degree 1 and $m-2$ unlabelled internal vertices of degree 3, and no other vertices. An *m-shape* is a tree having $m$ labelled external vertices of degree 1 and $m-2$ unlabelled internal vertices of degree 3, and no other vertices, *i.e.*, an $m$-shape is a labelling of an $m$-skeleton's external vertices by the labels $0, 1, \ldots, m-1$. When $m$ is clear from the context, we will refer to an $m$-shape simply as a shape. For notational convenience, we associate to each $m$-shape an arbitrary labelling of its $2m-3$ edges, with labels $1, \ldots, 2m-3$. This arbitrary choice of edge labelling is fixed once and for all. Thus an $m$-shape $\sigma$ is a labelling of an $m$-skeleton's external vertices together with a corresponding specification of edge labels. Let $\Sigma_m$ denote the set of $m$-shapes. There is a unique shape for $m = 2$ and $m = 3$, and $(2m - 5)!!$ distinct shapes for $m \geq 4$ (see [?, (5.96)] for a proof). In this notation, $(-1)!! = 1$ and $(2j+1)!! = (2j+1)(2j-1)!!$ for $j \geq 0$.

Let $m \geq 2$. Given a shape $\sigma \in \Sigma_m$, we associate to edge $j$ (oriented away from vertex 0) a nonnegative real number $t_j$ and a vector $y_j$ in $\mathbb{R}^d$. Writing $\vec{y} = (y_1, \ldots, y_{2m-3})$ and $\vec{t} = (t_1, \ldots, t_{2m-3})$, we define
$$a^{(m)}(\sigma; \vec{y}, \vec{t}) = \left(\sum_{i=1}^{2m-3} t_i\right) e^{-(\sum_{i=1}^{2m-3} t_i)^2/2} \prod_{i=1}^{2m-3} p_{t_i}(y_i) \qquad (2.4)$$

and
$$A^{(m)}(\sigma; \vec{y}) = \int_0^\infty dt_1 \cdots \int_0^\infty dt_{2m-3} \, a^{(m)}(\sigma; \vec{y}, \vec{t}). \qquad (2.5)$$

Then $\int_{\mathbb{R}^{d(2m-3)}} A^{(m)}(\sigma; \vec{y}) d\vec{y} = 1/(2m-5)!!$, so the sum of this integral over shapes $\sigma \in \Sigma_m$ is equal to 1. Let $\vec{k} \cdot \vec{y} = \sum_{j=1}^{2m-3} k_j \cdot y_j$, with each $k_j \in \mathbb{R}^d$. The Fourier integral transform $\hat{A}^{(m)}(\sigma; \vec{k}) = \int_{\mathbb{R}^{d(2m-3)}} A^{(m)}(\sigma; \vec{y}) e^{i\vec{k}\cdot\vec{y}} d\vec{y}$ is given by

$$\hat{A}^{(m)}(\sigma; \vec{k}) = \int_0^\infty dt_1 \cdots \int_0^\infty dt_{2m-3}\, \hat{a}^{(m)}(\sigma; \vec{k}, \vec{t}), \tag{2.6}$$

with

$$\hat{a}^{(m)}(\sigma; \vec{k}, \vec{t}) = \left(\sum_{i=1}^{2m-3} t_i\right) e^{-(\sum_{i=1}^{2m-3} t_i)^2/2} \prod_{i=1}^{2m-3} e^{-k_i^2 t_i/2}. \tag{2.7}$$

The $l^{\text{th}}$ moment measure $M^{(l)}$ for ISE can be written in terms of $A^{(l+1)}$, for $l \geq 1$. This is a deterministic measure which is absolutely continuous with respect to Lebesgue measure on $\mathbb{R}^{dl}$. The first moment measure $M^{(1)}$ has density $A^{(2)}(x)$. The second moment measure $M^{(2)}$ has density $\int A^{(3)}(y, x_1 - y, x_2 - y) d^d y$. In general, the density of $M^{(l)}$ at $x_1, \ldots, x_l$, for $l \geq 3$, is given by integrating $A^{(l+1)}(\sigma; \vec{y})$ over $\mathbb{R}^{d(l-1)}$ and then summing over the $(2l-3)!!$ shapes $\sigma$. Here $\vec{y}$ consists of integration variables $y_j$ corresponding to the edges $j$ on paths from vertex 0 to vertices of degree 3 in $\sigma$, and the other $y_a$ are fixed by the requirement that each external vertex $x_i$ is given by the sum of the $y_e$ over the edges $e$ connecting vertices 0 and $i$ in $\sigma$. Thus, the integration corresponds to integrating over the $l-1$ internal vertices, with the $l+1$ external vertices fixed at $0, x_1, \ldots, x_l$. For example, the contribution to the density of $M^{(3)}$ due to $\sigma_1$ of Figure **??** is $\int A^{(4)}(\sigma_1; y_1, x_1 - y_1, y_3, x_2 - y_1 - y_3, x_3 - y_1 - y_3) d^d y_1 d^d y_3$.

ISE and the functions (**??**) and (**??**) are further discussed in [**?**] (see also [**?**, **?**, **?**]). A construction of ISE as the scaling limit of branching random walk conditioned on the total size of the family tree, including a derivation of these functions, is given in [**?**].

# 3 Generating functions

For our applications to lattice trees and percolation, it will be essential to understand that the Fourier integral transforms of $a^{(m)}$ and $A^{(m)}$, $m = 2, 3, 4, \ldots$, occur in the asymptotic behaviour of certain generating function coefficients. This connection between ISE and generating functions was pointed out in [**?**].

The relevant generating functions are defined as follows. For $k \in \mathbb{R}^d$,

define $C^{(2)}_{z,\zeta}(k)$ and $c^{(2)}_{n,s}(k)$ by

$$C^{(2)}_{z,\zeta}(k) = \frac{2}{k^2 + 2^{3/2}\sqrt{1-z} + 2(1-\zeta)} = \sum_{s,n=0}^{\infty} c^{(2)}_{n,s}(k) z^n \zeta^s, \quad |\zeta|, |z| < 1, \tag{3.1}$$

where the square root has branch cut $[1,\infty)$ and is positive for real $z < 1$. For $m \geq 2$, given a shape $\sigma \in \Sigma_m$, to edge $j$ we associate $k_j \in \mathbb{R}^d$ and $\zeta_j \in \mathbb{C}$, with $|\zeta_j| < 1$. We write $\vec{k} = (k_1, \ldots, k_{2m-3})$ and $\vec{\zeta} = (\zeta_1, \ldots, \zeta_{2m-3})$, and define

$$C^{(m)}_{z,\vec{\zeta}}(\sigma; \vec{k}) = \prod_{j=1}^{2m-3} C^{(2)}_{z,\zeta_j}(k_j) = \sum_{s_1,\ldots,s_{2m-3}=0}^{\infty} \sum_{n=0}^{\infty} c^{(m)}_{n,\vec{s}}(\sigma; \vec{k}) z^n \prod_{j=1}^{2m-3} \zeta_j^{s_j}. \tag{3.2}$$

We write $b^{(m)}_{\vec{s}}(\sigma; \vec{k}) = \sum_{n=0}^{\infty} c^{(m)}_{n,\vec{s}}(\sigma; \vec{k})$ for the coefficient of $\prod_{j=1}^{2m-3} \zeta_j^{s_j}$ in $C^{(m)}_{1,\vec{\zeta}}(\sigma; \vec{k})$. Writing $\vec{1} = (1, \ldots, 1)$, we denote the coefficient of $z^n$ in $C^{(m)}_{z,\vec{1}}(\sigma; \vec{k})$ by $c^{(m)}_n(\sigma; \vec{k}) = \sum_{s_1,\ldots,s_{2m-3}=0}^{\infty} c^{(m)}_{n,\vec{s}}(\sigma; \vec{k})$.

The coefficients $b^{(m)}_{\vec{s}}(\sigma; \vec{k})$ are easily identified from the fact that $C^{(2)}_{1,\zeta}(k)$ is the sum of a geometric series in $\zeta$, namely

$$C^{(2)}_{1,\zeta}(k) = \frac{2}{k^2 + 2(1-\zeta)} = \sum_{s=0}^{\infty} \frac{1}{(1 + k^2/2)^{s+1}} \zeta^s. \tag{3.3}$$

Therefore $b^{(m)}_{\vec{s}}(\sigma; \vec{k}) = \prod_{j=1}^{2m-3}(1 + k_j^2/2)^{-(s_j+1)}$. For $t_j \in [0, \infty)$, the Fourier transform of the Brownian transition density (??) then emerges as the $m = 2$ case of the limit

$$\lim_{n \to \infty} b^{(m)}_{\lfloor \vec{t}n \rfloor}(\sigma; \vec{k} n^{-1/2}) = \prod_{j=1}^{2m-3} e^{-k_j^2 t_j / 2}. \tag{3.4}$$

Here $\lfloor \vec{t}n \rfloor$ denotes the vector with components $\lfloor t_j n \rfloor$.

For the ISE $m$-point function (??), we consider the generating function $C^{(m)}_{z,\vec{1}}(\sigma; \vec{k}) = \prod_{j=1}^{2m-3} 2(k_j^2 + 2^{3/2}\sqrt{1-z})^{-1}$. By Cauchy's theorem,

$$c^{(m)}_n(\sigma; \vec{k}) = \frac{1}{2\pi i} \oint_\Gamma C^{(m)}_{z,\vec{1}}(\sigma; \vec{k}) \frac{dz}{z^{n+1}}, \tag{3.5}$$

where $\Gamma$ is a circle centred at the origin with radius less than 1. By deforming the contour to the branch cut $[1, \infty)$ of the square root, it can be shown that for any $m \geq 2$, $k \in \mathbb{R}^d$,

$$c^{(m)}_n(\sigma; \vec{k} n^{-1/4}) \sim \frac{1}{\sqrt{2\pi}} n^{m-5/2} \hat{A}^{(m)}(\sigma; \vec{k}) \tag{3.6}$$

as $n \to \infty$. Here $f(n) \sim g(n)$ denotes $\lim_{n\to\infty} f(n)/g(n) = 1$. A proof of (??) is given in the proof of Theorem 1.1 of [?].

The functions $\hat{a}^{(m)}(\sigma; \vec{k}, \vec{t})$ arise from an appropriate joint limit of the coefficients $c_{n,\vec{s}}^{(m)}(\sigma; \vec{k})$. Namely, for $m \geq 2$,

$$c_{n,\lfloor \vec{t}n^{1/2} \rfloor}^{(m)}(\sigma; \vec{k}n^{-1/4}) \sim \frac{1}{\sqrt{2\pi}} \frac{1}{n} \hat{a}^{(m)}(\sigma; \vec{k}, \vec{t}) \qquad (3.7)$$

as $n \to \infty$. A proof is given in the proof of Theorem 1.2 of [?].

One might wonder at this point what any of this has to do with lattice trees or percolation. The connection is that some of these models' key thermodynamic functions have the form of the above generating functions in high dimensions, and this links them to ISE.

## 4 Lattice trees

A lattice tree in $\mathbb{Z}^d$ is a finite connected set of lattice bonds containing no cycles. For the nearest-neighbour model, the bonds are nearest-neighbour bonds $\{x, y\}$, $x, y \in \mathbb{Z}^d$, $\|x - y\|_1 = 1$. We will also consider "spread-out" lattice trees constructed from bonds $\{x, y\}$ with $0 < \|x - y\|_\infty \leq L$. The parameter $L$ will later be taken to be large but finite. We associate the uniform probability measure to the set of all $n$-bond lattice trees which contain the origin.

In this section, we will describe results showing that in high dimensions the scaling limit of lattice trees of size $n$, with space scaled by a multiple of $n^{-1/4}$, is ISE. Thus lattice trees in high dimensions behave like branching random walk.

We define the one-point function $t_n^{(1)}$ to be the number of $n$-bond lattice trees containing the origin, with $t_0^{(1)} = 1$. By a subadditivity argument, there is a positive constant $z_c$ (depending on $d$, and on $L$ for the spread-out model) such that $\lim_{n\to\infty} [t_n^{(1)}]^{1/n} = z_c^{-1}$.

Next, we would like to define the higher-point functions $t_n^{(m)}(\sigma; \vec{y}, \vec{s})$, for $m \geq 2$. These functions count lattice trees with a certain property. To describe this, we need some definitions. Let $\sigma \in \Sigma_m$, and associated to each edge $j$ in $\sigma$, let $y_j \in \mathbb{Z}^d$ and let $s_j$ be a nonnegative integer ($j = 1, \ldots, 2m - 3$). First, we introduce the notion of backbone. Given a lattice tree $T$ containing the sites $0, x_1, \ldots, x_{m-1}$, we define the *backbone* $B$ of $(T; 0, x_1, \ldots, x_{m-1})$ to be the subtree of $T$ spanning $0, x_1, \ldots, x_{m-1}$. There is an induced labelling of the external vertices of the backbone, in which vertex $x_l$ is labelled $l$. Ignoring vertices of degree 2 in $B$, this backbone is equivalent to a shape $\sigma_B$ or to its modification by contraction of one or more edges to a point. (In the latter case, as we will discuss further in Appendix **??**, the choice of $\sigma_B$ may not be unique.) Next, we need a notion of compatibility. Restoring vertices of degree 2 in $B$, let $b_j$ denote the length of the backbone path corresponding to edge $j$ of $\sigma_B$, with $b_j = 0$ for any contracted edge. We say that $(T; 0, x_1, \ldots, x_{m-1})$ is *compatible* with $(\sigma; \vec{y}, \vec{s})$ if $\sigma_B$ can be chosen (when not uniquely determined) such that $\sigma_B = \sigma$, if $b_j = s_j$ for all edges $j$ of $\sigma$, and if the backbone path corresponding to $j$ undergoes the displacement $y_j$ for all edges $j$ of $\sigma$.

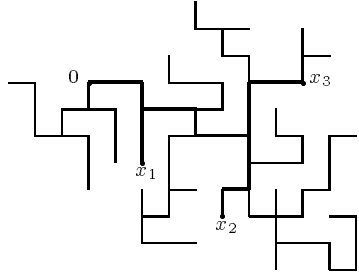

Figure 2: A 2-dimensional lattice tree contributing to $t_{78}^{(4)}(\sigma_1; \vec{y}, \vec{s})$, with $\sigma_1$ depicted in Figure **??**, $\vec{y} = ((2, -1), (0, -2), (4, -1), (-1, -3), (2, 2))$, $\vec{s} = (3, 2, 5, 4, 4)$.

Then we define $t_n^{(m)}(\sigma; \vec{y}, \vec{s})$ to be the number of $n$-bond lattice trees $T$, containing the origin, for which there are sites $x_1, \ldots, x_{m-1} \in T$ such that $(T; 0, x_1, \ldots, x_{m-1})$ is compatible with $(\sigma; \vec{y}, \vec{s})$. See Figure **??**. We also define
$$t_n^{(m)}(\sigma; \vec{y}) = \sum_{\vec{s}} t_n^{(m)}(\sigma; \vec{y}, \vec{s}), \tag{4.1}$$
where the sum over $\vec{s}$ denotes a sum over the nonnegative integers $s_j$. We will make use of Fourier transforms with respect to the $\vec{y}$ variables, for example,
$$\hat{t}_n^{(m)}(\sigma; \vec{k}) = \sum_{\vec{y}} t_n^{(m)}(\sigma; \vec{y}) e^{i\vec{k} \cdot \vec{y}}, \quad k_j \in [-\pi, \pi]^d. \tag{4.2}$$

For $m = 2, 3$ there is only one shape and we will sometimes omit it from the notation.

Define
$$G_{z,\vec{\zeta}}^{(m)}(\sigma; \vec{y}) = \sum_{n=0}^{\infty} \sum_{\vec{s}} t_n^{(m)}(\sigma; \vec{y}, \vec{s}) z^n \prod_{j=1}^{2m-3} \zeta_j^{s_j}. \tag{4.3}$$

The sum over $\vec{y} \in \mathbb{R}^{d(2m-3)}$ of (**??**) is finite for $|z| < z_c$ and $|\zeta_j| \le 1$, for all $m$.

In terms of critical exponents, the Fourier transform of the two-point function $G^{(2)}_{z,1}(y)$ is believed to behave asymptotically as

$$\hat{G}^{(2)}_{z_c,1}(k) \sim \frac{c_1}{k^{2-\eta}} \text{ as } k \to 0, \quad \hat{G}^{(2)}_{z,1}(0) \sim \frac{c_2}{(1-z/z_c)^\gamma} \text{ as } z \to z_c, \quad (4.4)$$

with the mean-field values $\eta = 0$ and $\gamma = \frac{1}{2}$ for $d > 8$. For $d > 8$, the simplest combination of the two asymptotic relations in (??) that could be hoped for is

$$\hat{G}^{(2)}_{z,1}(k) = \frac{C_1}{D_1^2 k^2 + 2^{3/2}(1-z/z_c)^{1/2}} + \text{error}, \quad (4.5)$$

where $C_1$ and $D_1$ are positive constants depending on $d$ and $L$. The error term is meant to be of lower order than the main term, in some suitable sense, as $k \to 0$ and $z \to z_c$.

An optimist expecting to find ISE and familiar with (??) and (??) could also hope that, for $d > 8$,

$$\hat{G}^{(2)}_{z,\zeta}(k) = \frac{C_1}{D_1^2 k^2 + 2^{3/2}(1-z/z_c)^{1/2} + 2T_1(1-\zeta)} + \text{error}, \quad (4.6)$$

and that there is an approximate independence of the form

$$\hat{G}^{(m)}_{z,\vec{\zeta}}(\sigma;\vec{k}) = v_1^{m-2} \prod_{j=1}^{2m-3} \hat{G}^{(2)}_{z,\zeta_j}(k_j) + \text{error}. \quad (4.7)$$

Here $v_1$ is a finite positive constant which translates the self-avoidance interactions of lattice trees into a renormalized vertex factor. For the nearest-neighbour model with $d$ sufficiently large, and for spread-out models for $d > 8$ with $L$ sufficiently large, relations of the form (??) and (??) are proved in [?], for all $m \geq 2$ if $\vec{\zeta} = \vec{1}$ and for $m = 2, 3$ for general $\vec{\zeta}$. The results given below arise as a consequence.

We define

$$p_n^{(m)}(\sigma; \vec{y}) = \frac{t_n^{(m)}(\sigma; \vec{y})}{\sum_{\sigma \in \Sigma_m} \hat{t}_n^{(m)}(\sigma; \vec{0})}, \quad (4.8)$$

which is a probability measure on $\Sigma_m \times \mathbb{Z}^{d(2m-3)}$. The following theorem, whose proof extends the methods of [?, ?], shows that (??) has the corresponding ISE density as its scaling limit in high dimensions. In its statement, the scaling of $\vec{k}$ by $D_1^{-1} n^{-1/4}$ corresponds to scaling down the lattice spacing by $D_1^{-1} n^{-1/4}$.

**Theorem 1** [?, ?] *Let $m \geq 2$ and $k_j \in \mathbb{R}^d$ ($j = 1, \ldots, 2m - 3$). For nearest-neighbour lattice trees in sufficiently high dimensions $d \geq d_0$, and for spread-out lattice trees with $d > 8$ and $L$ sufficiently large depending on $d$, there are constants $c_1, D_1$ depending on $d$ and $L$, such that*

$$\hat{t}_n^{(m)}(\sigma; \vec{k}D_1^{-1}n^{-1/4}) \sim c_1 n^{m-5/2} z_c^{-n} \hat{A}^{(m)}(\sigma; \vec{k}) \quad (n \to \infty). \tag{4.9}$$

*In particular,*
$$\lim_{n \to \infty} \hat{p}_n^{(m)}(\sigma; \vec{k}D_1^{-1}n^{-1/4}) = \hat{A}^{(m)}(\sigma; \vec{k}).$$

It is a corollary of Theorem **??** that high-dimensional lattice trees converge weakly to ISE, as we now explain. Given an $n$-bond lattice tree $T$ containing the origin, we define $\mu_n^T$ to be the probability measure on $\mathbb{R}^d$ which assigns mass $(n+1)^{-1}$ to the each of the $n+1$ points $xD_1^{-1}n^{-1/4}$, for $x \in T$. Let $M_1(\mathbb{R}^d)$ be the space of probability measures on $\mathbb{R}^d$. We then define a probability measure $\mu_n$ on $M_1(\mathbb{R}^d)$, supported on the $\mu_n^T$, by $\mu_n(\mu_n^T) = (t_n^{(1)})^{-1}$ for each $n$-bond $T$ containing 0. In this way, $n$-bond lattice trees induce a random probability measure on $\mathbb{R}^d$. Let $\dot{\mathbb{R}}^d$ denote the one-point compactification of $\mathbb{R}^d$, and let $M_1(\dot{\mathbb{R}}^d)$ denote the compact set of probability measures on $\dot{\mathbb{R}}^d$, under the topology of weak convergence. We regard $M_1(\mathbb{R}^d)$ as embedded in $M_1(\dot{\mathbb{R}}^d)$.

**Corollary 2** *For nearest-neighbour lattice trees in sufficiently high dimensions $d \geq d_0$, and for spread-out lattice trees with $d > 8$ and $L$ sufficiently large, $\mu_n$ converges weakly to $\mu_{\text{ISE}}$, as measures on $M_1(\dot{\mathbb{R}}^d)$.*

The weak convergence in Corollary **??** is the assertion that for any continuous function $F$ on $M_1(\dot{\mathbb{R}}^d)$,

$$\lim_{n \to \infty} \int_{M_1(\dot{\mathbb{R}}^d)} F(\nu) d\mu_n(\nu) = \int_{M_1(\dot{\mathbb{R}}^d)} F(\nu) d\mu_{\text{ISE}}(\nu). \tag{4.10}$$

The argument leading from Theorem **??** to Corollary **??** is presented in Appendix **??**.

For a more refined statement than Theorem **??**, we define

$$p_n^{(m)}(\sigma; \vec{y}, \vec{s}) = \frac{t_n^{(m)}(\sigma; \vec{y}, \vec{s})}{\sum_{\sigma \in \Sigma_m} \hat{t}_n^{(m)}(\sigma; \vec{0})}, \tag{4.11}$$

which is a probability measure on $\Sigma_m \times \mathbb{Z}^{d(2m-3)} \times \mathbb{Z}_+^{2m-3}$.

**Theorem 3** *[?, ?] Let $m = 2$ or $m = 3$, $k_j \in \mathbb{R}^d$, and $t_j \in (0, \infty)$ ($j = 1, \ldots, 2m - 3$). For nearest-neighbour lattice trees in sufficiently high dimensions $d \geq d_0$, and for spread-out lattice trees with $d > 8$ and $L$ sufficiently large depending on $d$, there is a constant $T_1$ depending on $d$ and $L$, such that*

$$\hat{t}_n^{(m)}(\sigma; \vec{k}D_1^{-1}n^{-1/4}, \lfloor \vec{t}T_1 n^{1/2} \rfloor) \sim c_1 T_1^{-(2m-3)} n^{-1} z_c^{-n} \hat{a}^{(m)}(\sigma; \vec{k}, \vec{t}) \quad (n \to \infty).$$

*In particular,*

$$\lim_{n \to \infty} (T_1 n^{1/2})^{2m-3} \hat{p}_n^{(m)}(\sigma; \vec{k}D_1^{-1}n^{-1/4}, \lfloor \vec{t}T_1 n^{1/2} \rfloor) = \hat{a}^{(m)}(\sigma; \vec{k}, \vec{t}). \quad (4.12)$$

We believe that Theorem ?? holds for all $m \geq 2$, but technical difficulties arise for $m \geq 4$ and the theorem has been proved only for $m = 2$ and $m = 3$. Theorem ?? indicates that, at least for $m = 2$ and $m = 3$, skeleton paths with length of order $n^{1/2}$ are typical. This is Brownian scaling, since distance is scaled as $n^{1/4}$. The statement of Theorem ?? for $m = 3$ in [11, 12] incorrectly included the case where $t_j = 0$ for one or two values of $j$, for which different constants occur, in fact, in the asymptotic formula for $\hat{t}_n^{(3)}(\sigma; \vec{k}D_1^{-1}n^{-1/4}, \lfloor \vec{t}T_1 n^{1/2} \rfloor)$.

We expect that the above results for lattice trees should apply also to lattice animals for $d > 8$, yielding ISE for their scaling limit for $d > 8$. This would be consistent with the general belief that lattice trees and lattice animals have the same scaling properties in all dimensions.

## 5  Percolation

Consider independent Bernoulli bond percolation on $\mathbb{Z}^d$, either nearest-neighbour or spread-out, with $p$ fixed and equal to its critical value $p_c$ [?]. Bonds are pairs $\{x, y\}$ of sites in $\mathbb{Z}^d$, with $\|x - y\|_1 = 1$ for the nearest-neighbour model and $0 < \|x - y\|_\infty \leq L$ for the spread-out model. Let $C(0)$ denote the random set of sites connected to 0, let $|C(0)|$ denote the cardinality of $C(0)$, and let

$$\tau^{(2)}(x; n) = P_{p_c}(C(0) \ni x, |C(0)| = n) \quad (5.1)$$

denote the probability at the critical point that the origin is connected to $x$ via a cluster containing $n$ sites. We define a generating function

$$\tau_z^{(2)}(x) = \sum_{n=1}^{\infty} \tau^{(2)}(x; n) z^n, \quad (5.2)$$

which converges absolutely if $|z| \leq 1$. Assuming no infinite cluster at $p_c$, $\tau_1^{(2)}(x)$ is the probability that 0 is connected to $x$.

The conventional definitions [?, Section 7.1] of the critical exponents $\eta$ and $\delta$ suggest that

$$\hat{\tau}_1^{(2)}(k) \sim \frac{c_1}{k^{2-\eta}} \text{ as } k \to 0, \quad \hat{\tau}_z^{(2)}(0) \sim \frac{c_2}{(1-z)^{1-1/\delta}} \text{ as } z \to 1, \quad (5.3)$$

but there is still no proof of existence of these exponents except in high dimensions. Assuming the mean-field values $\eta = 0$ and $\delta = 2$ above six dimensions, the simplest combination of the above asymptotic relations for $d > 6$ would be

$$\hat{\tau}_z^{(2)}(k) = \frac{C_2}{D_2^2 k^2 + 2^{3/2}(1-z)^{1/2}} + \text{error}, \quad (5.4)$$

for some constants $C_2$, $D_2$. This is analogous to (??). The following theorem shows that this behaviour is what does occur for sufficiently spread-out percolation above six dimensions.

**Theorem 4** *[?, ?] Let $k \in [-\pi, \pi]^d$, $z \in [0, 1)$. For spread-out percolation with $d > 6$ and $L$ sufficiently large, there are functions $\epsilon_1(z)$ and $\epsilon_2(k)$ with $\lim_{z \to 1} \epsilon_1(z) = \lim_{k \to 0} \epsilon_2(k) = 0$, and constants $C_2$ and $D_2$ depending on $d$ and $L$, such that*

$$\hat{\tau}_z^{(2)}(k) = \frac{C_2}{D_2^2 k^2 + 2^{3/2}(1-z)^{1/2}} [1 + \epsilon(z,k)] \quad (5.5)$$

*with $|\epsilon(z,k)| \leq \epsilon_1(z) + \epsilon_2(k)$.*

In view of (??), Theorem ?? is highly suggestive that ISE occurs as a scaling limit for percolation, but the control of the error term in (??) is too weak to obtain bounds on $\hat{\tau}^{(2)}(kD_2^{-1}n^{-1/4}; n)$ via contour integration. However, for the nearest-neighbour model in sufficiently high dimensions, better control of the error terms has been obtained, for *complex* $z$ with $|z| < 1$, leading to the following theorem. The theorem also gives a result for the three-point function

$$\tau^{(3)}(x, y; n) = P_{p_c}(x, y \in C(0), |C(0)| = n), \quad (5.6)$$

in terms of its Fourier transform

$$\hat{\tau}^{(3)}(k, l; n) = \sum_{x,y \in \mathbb{Z}^d} \tau^{(3)}(x, y; n) e^{ik \cdot x + il \cdot y}. \quad (5.7)$$

**Theorem 5** *[?, ?] Fix $k, l \in \mathbb{R}^d$ and any $\epsilon \in (0, \frac{1}{2})$. There is a $d_0$ such that for nearest-neighbour percolation with $d \geq d_0$, there are constants $C_2, D_2$ (depending on d) such that as $n \to \infty$*

$$\hat{\tau}^{(2)}(kD_2^{-1}n^{-1/4}; n) = \frac{C_2}{\sqrt{8\pi n}} \hat{A}^{(2)}(k)[1 + O(n^{-\epsilon})], \tag{5.8}$$

$$\hat{\tau}^{(3)}(kD_2^{-1}n^{-1/4}, lD_2^{-1}n^{-1/4}; n) = \frac{C_2}{\sqrt{8\pi}} n^{1/2} \hat{A}^{(3)}(k+l, k, l)[1 + O(n^{-\epsilon})]. \tag{5.9}$$

It follows from (??) that

$$P_{p_c}(|C(0)| = n) = n^{-1}\hat{\tau}^{(2)}(0; n) = C_2(8\pi)^{-1/2}n^{-3/2}[1 + O(n^{-\epsilon})]. \tag{5.10}$$

This shows that the critical exponent $\delta$, defined by $P_{p_c}(|C(0)| = n) \approx n^{-1-1/\delta}$, is given by $\delta = 2$ in high dimensions.

The variables in (??) are arranged schematically as:

```
            y
            |
            | l
  0 ________|________  x .
      k+l       k
```

To obtain (??), we work with the generating function

$$\hat{\tau}_z^{(3)}(k, l) = \sum_{n=1}^{\infty} \tau^{(3)}(k, l; n) z^n, \tag{5.11}$$

and prove that there is a positive constant $v_2$ such that

$$\hat{\tau}_z^{(3)}(k, l) = v_2 \hat{\tau}_z^{(2)}(k+l)\hat{\tau}_z^{(2)}(k)\hat{\tau}_z^{(2)}(l) + \text{error}. \tag{5.12}$$

An asymptotic relation in the spirit of (??), with $k = l = 0$, was conjectured for $d > 6$ already in [?].

We expect that Theorem ?? should extend to general $m$-point functions, for all $m \geq 2$, but this has not been proven. This is essentially the conjecture of [?] that the scaling limit of the incipient infinite cluster is ISE for $d > 6$. We now discuss this conjecture in more detail.

Given a site lattice animal $S$ containing $n$ sites, one of which is the origin, define the probability measure $\nu_n^S \in M_1(\mathbb{R}^d)$ to assign mass $n^{-1}$ to $xD_2^{-1}n^{-1/4}$, for each $x \in S$. We define $\nu_n$ to be the probability measure on $M_1(\mathbb{R}^d)$ which assigns probability $P_{p_c}(C(0) = S \mid |C(0)| = n)$ to $\nu_n^S$, for each $S$ as above. We regard the limit of $\nu_n$, as $n \to \infty$, as the scaling limit of the incipient infinite cluster. This is related to one of Kesten's definitions of the incipient infinite cluster [?], but here we are taking the lattice spacing to zero as $n \to \infty$. The conjecture of [?] is that, as in Corollary ?? above, $\nu_n$ converges weakly to $\mu_{\text{ISE}}$ for $d > 6$.

The conjecture is supported by Theorem ??. In fact, the characteristic functions $\hat{N}_n^{(1)}(k)$ and $\hat{N}_n^{(2)}(k,l)$ of the first and second moment measures $N_n^{(1)}$ and $N_n^{(2)}$ of $\nu_n$ are given by

$$\hat{N}_n^{(1)}(k) = \frac{\hat{\tau}^{(2)}(kD_2^{-1}n^{-1/4};n)}{\hat{\tau}^{(2)}(0;n)}, \tag{5.13}$$

$$\hat{N}_n^{(2)}(k,l) = \frac{\hat{\tau}^{(3)}(kD_2^{-1}n^{-1/4}, lD_2^{-1}n^{-1/4};n)}{\hat{\tau}^{(3)}(0,0;n)}, \tag{5.14}$$

and in high dimensions these converge respectively to the characteristic functions $\hat{A}^{(2)}(k)$ and $\hat{A}^{(3)}(k+l,k,l)$ of the corresponding ISE moments, by Theorem ??.

## 6  Oriented percolation

Consider independent oriented percolation on $\mathbb{Z}^d \times \mathbb{Z}_+$. Bonds are directed and are of the form $((x,n),(y,n+1))$, with $x,y \in \mathbb{Z}^d$ obeying $\|x-y\|_1 = 1$ for the nearest-neighbour model and obeying $0 < \|x-y\|_\infty \leq L$ for the spread-out model. Bonds are occupied with probability $p$. We write $(x,m) \to (y,n)$ if there is an oriented path from $(x,m)$ to $(y,n)$ consisting of occupied bonds, and define $C(x,m) = \{(y,n) : (x,m) \to (y,n)\}$. Let

$$\sigma^{(2)}((x,n);N) = P_{p_c}(C(0,0) \ni (x,n), |C(0,0)| = N) \tag{6.1}$$

denote the probability at the oriented percolation critical point that $(0,0)$ is connected to $(x,n)$ via a cluster containing $N$ sites. We denote the Fourier transform with respect to $x$ by

$$\hat{\sigma}^{(2)}((k,n);N) = \sum_{x \in \mathbb{Z}^d} \sigma^{(2)}((x,n);N)e^{ik\cdot x}, \quad k \in [-\pi,\pi]^d, \tag{6.2}$$

and define

$$\hat{\sigma}_{z,\zeta}^{(2)}(k) = \sum_{N=1}^{\infty} \sum_{n=0}^{\infty} \hat{\sigma}^{(2)}((k,n);N) z^N \zeta^n, \quad |z|,|\zeta| < 1. \tag{6.3}$$

The symmetry under $x \to -x$ is responsible for the absence of a term linear in $k$ in the denominators of (??) and (??). This symmetry applies also for oriented percolation, but there is no such symmetry for the time variable $n$ and a term linear in $(1-\zeta)$ should appear. Thus we expect that above the upper critical dimension, i.e., for $d+1 > 5$,

$$\hat{\sigma}^{(2)}_{z,\zeta}(k) = \frac{C_3}{D_3^2 k^2 + 2^{3/2}\sqrt{1-z} + 2T_3(1-\zeta)} + \text{error} \qquad (6.4)$$

as $(k, z, \zeta) \to (0, 1, 1)$. An upper bound for (??) of the form $(k^2 + |1-\zeta|)^{-1}$ was obtained for $z = 1$ in [?], for the nearest-neighbour model in sufficiently high dimensions and for sufficiently spread-out models when $d+1 > 5$. This is consistent with (??).

Apart from constants, the form of (??) is identical to the generating function $C'^{(2)}_{z,\zeta}(k)$ defined in (??). As in (??), if (??) accurately captures the behaviour of the two-point function, as $N \to \infty$ we would have

$$\hat{\sigma}^{(2)}((kD_3^{-1}N^{-1/4}, \lfloor tT_3 N^{1/2}\rfloor); N) \sim C_3 T_3^{-1} \frac{1}{\sqrt{8\pi N}} \hat{a}^{(2)}(k, t). \qquad (6.5)$$

This suggests ISE as the scaling limit, when time and space are scaled respectively by $N^{-1/2}$ and $N^{-1/4}$. The ISE time variable corresponds simply to the direction of orientation.

Consider now the limit in which the cluster size $N$ is summed over rather than fixed, with $n \to \infty$ and space scaled by $n^{-1/2}$. Summing over $N$ removes any conditioning on the cluster size, so SBM becomes relevant as the scaling limit, rather than ISE. According to the above picture, we can expect that

$$\hat{\sigma}^{(2)}_{1,\zeta}(k) = \frac{C_3}{D_3^2 k^2 + 2T_3(1-\zeta)} + \text{error.} \qquad (6.6)$$

As in (??), with sufficient control on the error (??) implies

$$2C_3^{-1}T_3 \lim_{n\to\infty} \sum_{N=1}^{\infty} \hat{\sigma}^{(2)}((kT_3^{1/2}D_3^{-1}n^{-1/2}, \lfloor tn \rfloor); N) = e^{-k^2 t/2}. \qquad (6.7)$$

In fact, (??)–(??) were proven in [?] for the nearest-neighbour model in sufficiently high dimensions and for sufficiently spread-out models when $d+1 > 5$. Work is in progress with Derbez and van der Hofstad to prove a corresponding result for higher-order connectivity functions, to obtain a stronger statement of convergence to SBM. This work in progress is based on the inductive method of [?], which bypasses the use of generating functions and the difficulties associated with their inversion.

The above picture relating SBM and oriented percolation can be contrasted with the results of [?] (see also [?]). In [?], it is shown that SBM arises as the scaling limit of a critical contact process for $d \geq 2$. The scaling limit of [?] is for the infinitely spread-out contact process, in the limit $L \to \infty$ (sometimes called the Kac limit). This is a mean-field limit, for which the non-gaussian behaviour expected below $d + 1 = 5$ when $L$ is finite is no longer relevant.

# 7 The lace expansion

The method of proof of the above results is based on the lace expansion, which was first introduced in [?] in the context of self-avoiding walks. Reviews of work on the lace expansion prior to the work described in this paper can be found in [?, ?]. The extensions required to prove the results of Sections ?? and ?? make use of a double lace expansion and it is beyond the scope of this paper to indicate any details. Details can be found in [?, ?, ?].

# A  Proof of Corollary ??

In this appendix, we show how Corollary ?? follows from Theorem ??. The corollary follows in a straightforward way via [?, Lemma 2.4.1(b)], which asserts that weak convergence of moment measures implies weak convergence of random probability measures (on a compact set). However, there is one subtlety. This point was overlooked in [?, ?], and we take this opportunity to clarify it.

For $l \geq 1$, let $s_n^{(l+1)}(x_1, \ldots, x_l)$ denote the number of $n$-bond lattice trees containing the lattice sites $0, x_1, \ldots, x_l$. To abbreviate the notation, we will write $\tilde{x} = (x_1, \ldots, x_l)$. The $l^{\text{th}}$ moment measure $M_n^{(l)}$ of $\mu_n$ is the deterministic probability measure on $\mathbb{R}^{dl}$ which places mass

$$r_n^{(l+1)}(\tilde{x}) = \frac{1}{(n+1)^l} \frac{1}{t_n^{(1)}} s_n^{(l+1)}(\tilde{x}) \tag{A.1}$$

at $\tilde{x} D_2^{-1} n^{-1/4}$, for $\tilde{x} \in \mathbb{Z}^{dl}$. The characteristic function $\hat{M}_n^{(l)}(k)$ of $M_n^{(l)}$ is given by

$$\hat{M}_n^{(l)}(\tilde{k}) = \hat{r}_n^{(l+1)}(\tilde{k} D_2^{-1} n^{-1/4}), \tag{A.2}$$

where, writing $\tilde{k} = (k_1, \ldots, k_l)$ and $\tilde{k} \cdot \tilde{x} = k_1 \cdot x_1 + \cdots + k_l \cdot x_l$,

$$\hat{r}_n^{(l+1)}(\tilde{k}) = \sum_{\tilde{x}} r_n^{(l+1)}(\tilde{x}) e^{i \tilde{k} \cdot \tilde{x}}. \tag{A.3}$$

Since $\hat{s}_n^{(l+1)}(\tilde{0}) = (n+1)^l t_n^{(1)}$, we have

$$\hat{M}_n^{(l)}(\tilde{k}) = \frac{\hat{s}_n^{(l+1)}(\tilde{k}D_1^{-1}n^{-1/4})}{\hat{s}_n^{(l+1)}(\tilde{0})}. \quad (A.4)$$

To prove convergence of the moment measures of $\mu_n$ to those of ISE, it suffices to show that, for each $l \geq 1$, $\hat{M}_n^{(l)}(k)$ converges to the characteristic function $\hat{M}^{(l)}(k)$ of the corresponding ISE moment measure described under (??). For $l = 1$, this is an immediate consequence of (??) and Theorem ??, since $\hat{s}_n^{(2)}(k) = \hat{t}_n^{(2)}(k)$ and $\hat{M}^{(1)}(k) = \hat{A}^{(2)}(k)$. Similarly, for $l = 2$, there is a unique shape and $\hat{s}_n^{(3)}(k_1, k_2) = \hat{t}_n^{(3)}(k_1 + k_2, k_1, k_2)$. Since $\hat{M}^{(2)}(k_1, k_2) = \int A^{(3)}(y, x_1 - y, x_2 - y)e^{ik_1 \cdot x_1}e^{ik_2 \cdot x_2}d^d y d^d x_1 d^d x_2 = \hat{A}^{(3)}(k_1 + k_2, k_1, k_2)$, convergence of the second moments follows directly from Theorem ??.

The convergence of the third and higher moments follows similarly, apart from one detail. For $l \geq 3$, there is more than one shape, and

$$\hat{M}^{(l)}(\tilde{k}) = \sum_{\sigma \in \Sigma_{l+1}} \hat{A}^{(l+1)}(\sigma; \vec{k}) \quad (A.5)$$

with each of the $2l - 1$ components of $\vec{k}$ given by a specific linear combination (depending on $\sigma$) of the $l$ components of $\tilde{k}$. For example, for $l = 3$ and the shape $\sigma_1$ of Figure ??, $(\sigma_1; \vec{k}) = (\sigma_1; k_1 + k_2 + k_3, k_1, k_2 + k_3, k_2, k_3)$. If it were the case that $\hat{s}_n^{(l+1)}(\tilde{k})$ were equal to $\sum_{\sigma \in \Sigma_{l+1}} \hat{t}_n^{(l+1)}(\sigma; \vec{k})$, convergence of all moments would be immediate since Theorem ?? implies that

$$\lim_{n \to \infty} \frac{\sum_{\sigma \in \Sigma_{l+1}} \hat{t}_n^{(l+1)}(\sigma; \vec{k}D_1^{-1}n^{-1/4})}{\sum_{\sigma \in \Sigma_{l+1}} \hat{t}_n^{(l+1)}(\sigma; \vec{0})} = \sum_{\sigma \in \Sigma_{l+1}} \hat{A}^{(l+1)}(\sigma; \vec{k}). \quad (A.6)$$

But $\hat{s}_n^{(l+1)}(\tilde{k})$ is not equal to $\sum_{\sigma \in \Sigma_{l+1}} \hat{t}_n^{(l+1)}(\sigma; \vec{k})$, because it is not the case that $s_n^{(l+1)}(\tilde{x})$ is equal to the sum of $t_n^{(l+1)}(\sigma; \vec{y})$ over all $(\sigma; \vec{y})$ that are consistent with $\tilde{x}$ in the sense that the $x_i$ are given by the sum of the $y_j$ as prescribed by the shape $\sigma$. The discrepancy arises from degenerate lattice tree configurations, containing sites $x_1, \ldots, x_l$, which can correspond to more than one choice of $(\sigma; \vec{y})$. These configurations can only occur when $l \geq 3$ and at least one $y_j$ is zero.

For example, there is a unique 1-bond lattice tree containing 0 and the site $e_1 = (1, 0, \ldots, 0)$, and hence $s_1^{(4)}(0, 0, e_1) = 1$. However, this lattice tree containing the sites $x_1 = x_2 = 0$, $x_3 = e_1$ contributes to each of $t_1^{(4)}(\sigma_1; 0, 0, 0, 0, e_1)$, $t_1^{(4)}(\sigma_2; 0, 0, 0, 0, e_1)$ and $t_1^{(4)}(\sigma_3; 0, e_1, 0, 0, 0)$. See Figure ??. Thus it is not the case, in general, that $s_n^{(l+1)}(\tilde{x})$ is given by the sum of $t_n^{(l+1)}(\sigma;\vec{y})$ over all corresponding $(\sigma;\vec{y})$. The assertion of [?, (3.4)] and [?, (1.11)] that $\sum_{\sigma \in \Sigma_{l+1}} \hat{t}_n^{(l+1)}(\sigma;\vec{0})$ equals $(n+1)^l t_n^{(1)}$ implicitly assumed uniqueness of $(\sigma;\vec{y})$ and is incorrect for $l \geq 3$. This false assertion was not needed in [?, ?], as it can be replaced by [?, (2.14)-(2.15)] with $\vec{k} = \vec{0}$ (i.e., (??) above) and [?, (1.12)] to conclude that

$$\sum_{\sigma \in \Sigma_{l+1}} \hat{t}_n^{(l+1)}(\sigma;\vec{0}) \sim c_1 n^{l-3/2} z_c^{-n} \sim (n+1)^l t_n^{(1)}, \tag{A.7}$$

which is sufficient for [?, ?]. The degenerate cases appear in error terms to (??) and do not affect the leading behaviour.

In view of (??)–(??), to prove convergence of the $l^{\text{th}}$ moments, for $l \geq 3$, it suffices to show that

$$\left| \hat{s}_n^{(l+1)}(\tilde{k}) - \sum_{\sigma \in \Sigma_{l+1}} \hat{t}_n^{(l+1)}(\sigma;\vec{k}) \right| \leq O(n^{l-2} z_c^{-n}). \tag{A.8}$$

This difference then constitutes an error term, down by $n^{-1/2}$ compared to $\hat{s}_n^{(l+1)}(\tilde{k})$, by Theorem ??. The remainder of the proof is devoted to obtaining (??).

Let $l \geq 3$, and recall the definition of compatibility above (??). If the backbone of $(T;0,x_1,\ldots,x_l)$ comprises $2l-1$ *nontrivial* paths (each having length greater than zero), then $\tilde{x}$ induces a labelling of the external vertices of an $(l+1)$-skeleton and there is therefore a unique compatible $(\sigma;\vec{y},\vec{s})$. Whether or not the backbone comprises $2l-1$ nontrivial paths, given $(\sigma;\vec{y},\vec{s})$ compatible with the backbone, the $2l-1$ backbone displacements $\vec{y}$ and their lengths $\vec{s}$ (possibly zero) are uniquely determined by $\sigma$ and $(T;0,x_1,\ldots,x_l)$. Nonuniqueness of $(\sigma;\vec{y},\vec{s})$ thus requires at least one of the backbone paths to be trivial, and, in such a degenerate case, the maximum possible number of compatible choices for $(\sigma;\vec{y},\vec{s})$ is the number of shapes, which is $(2l-3)!!$. Let $u_n^{(l+1)}(\tilde{x})$ denote the number of $n$-bond lattice trees for which each of the $2l-1$ backbone paths is nontrivial, and let $e_n^{(l+1)}(\tilde{x})$ denote the number of $n$-bond lattice trees for which at least one backbone path has a zero displacement. Then $s_n^{(l+1)}(\tilde{x}) = u_n^{(l+1)}(\tilde{x}) + e_n^{(l+1)}(\tilde{x})$, and, for $l \geq 3$,

$$\left| \hat{s}_n^{(l+1)}(\tilde{k}) - \sum_{\sigma \in \Sigma_{l+1}} \hat{t}_n^{(l+1)}(\sigma;\vec{k}) \right| \leq [(2l-3)!! - 1]\hat{e}_n^{(l+1)}(\tilde{0}). \tag{A.9}$$

It suffices to argue that the right side of (??) is at most $O(n^{l-2} z_c^{-n})$. For this, we introduce the generating function $E^{(l+1)}(z) = \sum_n \hat{e}_n^{(l+1)}(\tilde{0}) z^n$. Let $\chi(z) = \sum_x G_z^{(2)}(x)$. It can be shown using standard bounds that $|E^{(l+1)}(z)| \leq O(\chi(|z|)^{2l-2})$, where the power $2l-2$ arises because at least one of the $2l-1$ backbone paths is trivial. Using the methods of [?], this can be refined to $|E^{(l+1)}(z)| \leq O(|\chi(z)|^3 \chi(|z|)^{2l-5})$, uniform in $|z| < z_c$. It follows from [?, (1.12)] that $|E^{(l+1)}(z)| \leq O(|1-z/z_c|^{-3/2}(1-|z|/z_c)^{-l+5/2})$. Then [?, Lemma 3.2(i)] implies the desired bound $\hat{e}_n^{(l+1)}(\tilde{0}) \leq O(n^{l-2} z_c^{-n})$. □

## Acknowledgements

This work was supported in part by NSERC. It is a pleasure to thank Eric Derbez and Takashi Hara for the enjoyable collaborations that led to the results described in this paper, and Christian Borgs, Jennifer Chayes, Remco van der Hofstad and Ed Perkins for valuable conversations. This paper was written primarily during a visit to Microsoft Research.

Department of Mathematics and Statistics, McMaster University, Hamilton, ON, Canada L8S 4K1
slade@math.mcmaster.ca

Address after July 1, 1999:
Department of Mathematics, University of British Columbia, Vancouver, BC, Canada V6T 1Z2
slade@math.ubc.ca